\documentclass[11pt]{article}
\usepackage{amsmath,amssymb,theorem}
\usepackage{graphicx}
\usepackage{subfigure}
\usepackage{psfrag}
\newtheorem{thm}{Theorem}[section]
\newtheorem{conj}[thm]{Conjecture}
\newtheorem{cor}[thm]{Corollary}
\newtheorem{lem}[thm]{Lemma}
\theorembodyfont{\rmfamily}
\def\pf{\bigskip\noindent {\bf Proof.}~~}

\def\dfn#1{{\sl #1}}

\def\less{\backslash}

\def\pf{\bigskip\noindent {\emph{Proof.}}~~}

\newcounter{counter}
\def\pr#1{(\ref{#1})}

\begin{document}
\title{Erd\H{o}s-Lov\'asz Tihany Conjecture  for graphs with forbidden holes}
\author{
 Zi-Xia Song\thanks{E-mail address: Zixia.Song@ucf.edu}
\\
Department  of Mathematics\\
 University of Central Florida\\
Orlando, FL 32816, USA\\
}
\maketitle
\begin{abstract}

  A \dfn{hole}  in a graph   is an induced cycle of length at least $4$.  Let $s\ge2$ and $t\ge2$ be integers.  A graph $G$ is \dfn{$(s,t)$-splittable} if $V(G)$ can be partitioned into two sets $S$ and $T$ such that $\chi(G[S ]) \ge s$ and $\chi(G[T ]) \ge t$. 
  The well-known Erd\H{o}s-Lov\'asz Tihany Conjecture from 1968 states that every  graph $G$ with $\omega(G) < \chi(G) = s + t - 1$ is $(s,t)$-splittable. 
 This conjecture is hard, and  few related results are known.  However, it has been verified to be true for line graphs, quasi-line graphs, and graphs with independence number $2$.    In this paper, we establish more evidence for  the Erd\H{o}s-Lov\'asz Tihany Conjecture by showing that every graph  $G$ with $\alpha(G)\ge3$,  $\omega(G) < \chi(G) = s + t - 1$,  and    no hole of length between $4$ and $2\alpha(G)-1$ is $(s,t)$-splittable, where $\alpha(G)$ denotes the independence number of a graph $G$.

\end{abstract}

\section{Introduction}
\baselineskip 16pt

All graphs considered in this paper are finite and without loops or multiple edges.
For a graph $G$, we   use $V(G)$ to denote the vertex set, $E(G)$ the edge set, $|G|$ the number of vertices,   $\Delta(G)$  the maximum degree, $\alpha(G)$ the independence number, $\omega(G)$ the clique number, $\chi(G)$ the chromatic number, and $\overline{G}$ the complement of $G$.
For a vertex $x \in V(G)$, we will use $N(x)$ to denote the set of vertices in $G$ which are adjacent to $x$.
We define $N[x] = N(x) \cup \{x\}$ and $d(x) = |N(x)|$.
Given vertex sets $A, B \subseteq V(G)$, we say that $A$ is \dfn{complete to} (resp. \dfn{anti-complete to}) $B$ if for every $a \in A$ and every $b \in B$, $ab \in E(G)$ (resp. $ab \notin E(G)$).
The subgraph of $G$ induced by $A$, denoted $G[A]$, is the graph with vertex set $A$ and edge set $\{xy \in E(G) : x, y \in A\}$. We denote by $B \less A$ the set $B - A$,   and $G \less A$ the subgraph of $G$ induced on $V(G) \less A$, respectively.
If $A = \{a\}$, we simply write $B \less a$    and $G \less a$, respectively.  A graph $H$ is an \dfn{induced subgraph} of a graph $G$ if $V(H)\subseteq V(G)$   and $H=G[V(H)]$.  A graph $G$ is \dfn{claw-free} if $G$ does not contain $K_{1,3}$ as an induced subgraph. 
Given two isomorphic graphs $G$ and $H$, we may (with a slight but common abuse of notation) write $G = H$.    A cycle with $t \ge 3$ vertices is denoted by $C_t$.    We use the convention   ``A :="  to mean that $A$ is defined to be
the right-hand side of the relation.\medskip

Let $s\ge2$ and $t\ge2$ be integers. A graph $G$ is \dfn{$(s,t)$-splittable} if $V(G)$ can be partitioned into two sets $S$ and $T$ such that $\chi(G[S ]) \ge s$ and $\chi(G[T ]) \ge t$.  In 1966, the following conjecture of Lov\'asz was published by Erd\H os~\cite{Erdos1968} and is known as the Erd\H os-Lov\'asz Tihany Conjecture.

\begin{conj}\label{conj:ELTC}
Let $G$ be a   graph with $\omega(G)<\chi(G)=s+t-1$, where   $s\ge2$ and $t\ge2$ are integers.  Then $G$ is $(s,t)$-splittable.  
\end{conj}

To date, Conjecture~\ref{conj:ELTC} has been shown to be true only for values of $(s, t) \in \{ (2, 2), \allowbreak (2, 3), \allowbreak (2, 4), \allowbreak (3, 3), \allowbreak (3, 4), \allowbreak (3, 5) \}$.
The case $(2, 2)$ is trivial.
The case $(3, 3)$ was shown by Brown and Jung in 1969~\cite{Brown1969}.
Mozhan~\cite{Mozhan1987} and Stiebitz~\cite{Stiebitz1987a} each independently showed the case $(2, 4)$ in 1987.
The cases $(2, 3)$,   $(3, 4)$, and $(3, 5)$ were also settled by Stiebitz in 1987 \cite{Stiebitz1987b}.
Recent work on the Erd\H os-Lov\'asz Tihany Conjecture has focused on proving the conjecture for certain classes of graphs. A graph $G$  is a \dfn{quasi-line graph} if for every vertex $v\in V(G)$, the  set of neighbors of $v$ can be covered by two  cliques, namely the vertex set of the neighborhood of $v$ can be partitioned into two cliques.  By
definition, quasi-line graphs are claw-free. Recently, quasi-line graphs attracted more attention (see \cite{CF2007, CF2008, CP}). In particular,
Chudnovsky and Seymour \cite{CP} gave a constructive characterization of quasi-line graphs.
Kostochka and Stiebitz~\cite{Kostochka2008} showed that Conjecture~\ref{conj:ELTC} holds for line graphs.
Balogh, Kostochka, Prince, and Stiebitz~\cite{Balogh2009} then showed that Conjecture~\ref{conj:ELTC} holds for all quasi-line graphs,  and all graphs $G$ with $\alpha(G) = 2$.  
\medskip


\begin{thm}[\cite{Balogh2009}]\label{quasiline}
 Let $G$ be a graph with  $\omega(G)<\chi(G)=s+t-1$, where   $s\ge2$ and $t\ge2$ are integers. If $G$ is a quasi-line graph or $\alpha(G)=2$, then   $G$ is $(s,t)$-splittable.
\end{thm}
More recently, Chudnovsky, Fradkin, and Plumettaz~\cite{Chudnovsky2013} proved the following slight weaking of Conjecture~\ref{conj:ELTC} for claw-free graphs, the proof of which is long and relies heavily on the structure theorem for claw-free graphs developed by Chudnovsky and Seymour~\cite{Chudnovsky2008}.

\begin{thm}\label{thm:ELTCClawFree}
Let $G$ be a claw-free graph with $\chi(G) > \omega(G)$.
Then there exists a clique $K$ with $|V(K)| \le 5$ such that $\chi(G \setminus V(K)) > \chi(G) - |V(K)|$.
\end{thm}

The most recent result related to the Erd\H os-Lov\'asz Tihany Conjecture is due to Stiebitz~\cite{Stiebitz2016}, who showed that for integers $s, t \ge 2$, any graph $G$ with $\omega(G) < \chi(G) = s + t - 1$ contains disjoint subgraphs $G_1$ and $G_2$ of $G$ with either $\chi(G_1) \ge s$ and $\text{col}(G_2) \ge t$, or $\text{col}(G_1) \ge s$ and $\chi(G_2) \ge t$, where $\text{col}(H)$ denotes the coloring number of a graph $H$.
\medskip

If we restrict $s = 2$ in Conjecture~\ref{conj:ELTC}, then the Erd\H os-Lov\'asz Tihany Conjecture states that for any graph $G$ with $\chi(G) > \omega(G) \ge 2$, there exists an edge $xy \in E(G)$ such that $\chi(G \less \{x, y\}) \ge \chi(G) - 1$.
To prove this special case of Conjecture~\ref{conj:ELTC}, suppose for a contradiction that no such an edge exists.
Then $\chi(G \less \{x, y\}) = \chi(G) - 2$ for every edge $xy \in E(G)$.
This motivates the definition of double-critical graphs.
A connected graph $G$ is \dfn{double-critical} if for every edge $xy \in E(G)$, $\chi(G \less \{x, y\}) = \chi(G) - 2$.
A graph $G$ is \dfn{$t$-chromatic} if $\chi(G) = t$.
We are now ready to state the following conjecture, which is referred to as the \dfn{Double-Critical Graph Conjecture}, due to Erd\H os and Lov\'asz~\cite{Erdos1968}.

\begin{conj}\label{conj:DC}
Let $G$ be a double-critical, $t$-chromatic graph.
Then $G = K_t$.
\end{conj}

Since Conjecture~\ref{conj:DC} is a special case of Conjecture~\ref{conj:ELTC}, it has been settled in the affirmative for $t \le 5$~\cite{Mozhan1987, Stiebitz1987a}, for line graphs~\cite{Kostochka2008}, for quasi-line graphs,  and for  graphs with independence number two~\cite{Balogh2009}.
Representing a weakening of Conjecture~\ref{conj:DC}, Kawarabayashi, Pedersen, and Toft~\cite{Kawarabayashi2010} have shown that any double-critical, $t$-chromatic graph contains $K_t$ as a minor for $t \in \{6, 7\}$.
As a further weakening, Pedersen~\cite{Pedersen2011} showed that any double-critical, $8$-chromatic graph contains $K_8^-$ as a minor.
Albar and Gon\c calves~\cite{Albar2015} later proved that any double-critical, 8-chromatic graph contains $K_8$ as a minor.
Their proof is computer-assisted.
Rolek and the present author~\cite{RolekSong2018} gave a computer-free proof of the same result and further showed that any double-critical, $t$-chromatic graph contains $K_9$ as a minor for all $t \ge 9$.
We note here that Theorem~\ref{thm:ELTCClawFree} does not completely settle Conjecture~\ref{conj:DC} for all claw-free graphs.
Recently, Huang and Yu~\cite{Huang2016} proved that the only double-critical, $6$-chromatic, claw-free graph is $K_6$. Rolek and the present author~\cite{RolekSong2017} further proved that the only double-critical, $t$-chromatic, claw-free graph is $K_t$ for all  $t\le8$. \medskip

  In this paper, we establish more evidence for the  Erd\H{o}s-Lov\'asz Tihany Conjecture.  By Theorem~\ref{quasiline},  Erd\H{o}s-Lov\'asz Tihany Conjecture holds for graphs $G$ with  $\alpha(G)=2$ but remains unknown for graphs $G$  with  $\alpha(G)\ge3$.  Let $\mathcal{F}$ be a family of graphs. A graph is \dfn{$\mathcal{F}$-free} if it does not contain any $F\in\mathcal{F}$ as an induced subgraph. We prove the following main result. 

\begin{thm}\label{main}
 Let $G$ be a   graph with $\alpha(G)\ge3$ and  $\omega(G)<\chi(G)=s+t-1$, where   $s\ge2$ and $t\ge2$ are integers.    If $G$ is $\{C_4, \, C_5, C_6, \dots, C_{2\alpha(G)-1}\}$-free, then $G$ is $(s,t)$-splittable. 
\end{thm}

We prove Theorem~\ref{main} in Section~\ref{mainresult}. Our proof of Theorem~\ref{main} relies on Theorem~\ref{quasiline} and  the following  well-known Strong Perfect Graph Theorem~\cite{spgt}.    A \dfn{hole} in a graph $G$ is an induced cycle of length at least $4$.  An \dfn{antihole} in $G$ is an induced subgraph isomorphic to the  complement of  a hole. A graph $G$ is \dfn{perfect} if $\chi(H)=\omega(H)$ for every induced subgraph $H$ of $G$.

\begin{thm}[\cite{spgt}]\label{spgt}
A graph  is perfect if and only if it has no odd hole and no odd antihole.
\end{thm}

We shall need the following corollary which was   observed   in \cite{ThomasSong2017}.  
\begin{cor}[\cite{ThomasSong2017}]\label{spgt1}
If $G$  is  $\{C_4, C_5, C_7, \dots,  C_{2\alpha(G)+1}\}$-free, then $G$ is perfect.
\end{cor}

\section{Proof of Theorem~\ref{main}}\label{mainresult}

The main idea in the proof of Theorem~\ref{main} is similar to that used in the proof of a result  of Thomas and the present author (see Theorem 2.3 in \cite{ThomasSong2017}), which states that Hadwiger's Conjecture~\cite{had} is true for $\{C_4, C_5, C_6, \dots, C_{2\alpha(G)-1}\}$-free  graphs $G$ with $\alpha(G)\ge3$.   We will show that any minimal counterexample to Theorem~\ref{main} is a quasi-line graph. First, we establish Lemma~\ref{quasi-line}, noting that it is deduced  from the proof of Theorem 2.3 given in \cite{ThomasSong2017}.  Our hope is that if a conjecture has been proven true for quasi-line graphs, then Lemma~\ref{quasi-line} might be used to demonstrate that such a conjecture also holds for $\{C_4, C_5, C_6, \dots, C_{2\alpha(G)-1}\}$-free graphs $G$ with $\alpha(G)\ge3$.

\begin{lem}\label{quasi-line}
Let $G$ be a $\{C_4, C_5, C_6, \dots, C_{2\alpha(G)-1}\}$-free graph  with  $\alpha(G)\ge3$ and $\Delta(G)\le |G|-2$. If $G$   contains  an induced cycle of length $2\alpha(G)+1$, then $G$ is a quasi-line graph. 
\end{lem}

\pf  Let $G$ be given as in the statement.  Let $C$ be an induced cycle of length $2\alpha+1$ in $G$  with vertices  $ v_0, v_1, \dots, v_{2\alpha}$ in order, where $\alpha :=\alpha(G)$.    We next prove a series of claims.  \medskip

\noindent \refstepcounter{counter}\label{e:Lnbr} (\arabic{counter})  For every  $w\in V(G\less C)$,  either $w$ is complete to $C$, or $w$ is adjacent to exactly three consecutive vertices on $C$, or $w$ is adjacent to exactly four consecutive vertices on $C$. 

\pf Since  $\alpha(G)=\alpha$, we see that $w$ is adjacent to at least one vertex on $C$. Suppose that $w$ is not complete to $C$. We may assume that  $wv_0\notin E(G)$ but $wv_1\in E(G)$. Then $w$ is not adjacent to $v_{2\alpha}, v_{2\alpha -1}, \dots , v_5$ because  $G$ is $\{C_4, \,C_5, \, \dots , \, C_{2\alpha-1}\}$-free.    If $wv_4\in E(G)$, then $w$ must be adjacent to $v_2, v_3$ because  $G$ is $\{C_4, \, C_5\}$-free. If $wv_4\notin E(G)$, then again $w$ must be adjacent to $v_2, v_3$  because $\alpha(G)=\alpha$.  Thus $w$ is adjacent to either  $v_1, v_2, v_3$  or $v_1, v_2,v_3, v_4$  on $C$,  as desired. \hfill\vrule height3pt width6pt depth2pt\medskip

Let $J$  denote the set of vertices in $G$ that are complete to $C$.  For each $i\in I : =\{0,1,\dots, 2\alpha\}$, let  $A_i\subseteq V(G\less C)$ (possibly empty) denote the set of vertices in $G$ adjacent to precisely $v_i, v_{i+1}, v_{i+2}$ on $C$, and let $B_i\subseteq V(G\less C)$ (possibly empty) denote the set of vertices in $G$ adjacent to precisely $v_i, v_{i+1}, v_{i+2}, v_{i+3}$ on $C$, where all arithmetic on indices here and henceforth is done modulo $2\alpha+1$. By \pr{e:Lnbr}, $\{J, V(C), A_0, A_1, \dots, A_{2\alpha}, B_0, B_1,  \dots, B_{2\alpha}\}$ partitions $V(G)$. \\

\noindent \refstepcounter{counter}\label{e:LJ} (\arabic{counter})
   $J=\emptyset$.

 \pf Suppose that $J\not=\emptyset$.  Then $G[J]$ must be a clique because   $G$ is $C_4$-free. Let   $a\in J$.  Since $\Delta(G)\le |G|-2$, there must exist  a vertex $b\in V(G)\less (V(C)\cup J)$  such that $ab\notin E(G)$. By \pr{e:Lnbr}, we may assume that $b$ is adjacent to $v_0, v_1, v_2$. But then $G[\{a, v_0, b, v_2\}]$ is an induced $C_4$ in $G$, a contradiction. \hfill\vrule height3pt width6pt depth2pt\medskip

  The fact that $\alpha(G) = \alpha$ implies that \medskip

\noindent \refstepcounter{counter}\label{e:ABclique} (\arabic{counter}) For each $i\in I$, both $G[A_i]$ and $G[B_i]$ are cliques; and $A_i$ is complete to $A_{i-1}\cup A_{i+1}$.\medskip

Since $G$ is $\{C_4, \,C_5, \, \dots , \, C_{2\alpha-1}\}$-free, one can easily check that \\

\noindent \refstepcounter{counter}\label{e:LAac} (\arabic{counter})
For each $i\in I$, $A_i$ is  anti-complete to each $A_j$, where $j\in I\less\{i-2, i-1, i, i+1, i+2\}$; and \medskip

\noindent \refstepcounter{counter}\label{e:LBac} (\arabic{counter})
For each $i\in I$, $B_i$ is complete to $B_{i-1}\cup A_{i}\cup A_{i+1}\cup B_{i+1}$ and anti-complete to each $B_j$, where $j\in I\less\{ i-1, i, i+1\}$.\medskip

We shall also  need the following:\medskip

\noindent \refstepcounter{counter}\label{e:Loneneighbor} (\arabic{counter})
For all $i\in I$, if $B_i\not=\emptyset$, then  $B_{j}=\emptyset$ for any $j\in I\less\{i-2, i-1, i, i+1, i+2\}$.

\pf    Suppose $B_i\not=\emptyset$ and $B_{j}\not=\emptyset$ for some $j\in I\less\{ i+2, i+1, i, i-1, i-2\}$.  We may assume that $j>i$. Let $a\in B_i$ and $b\in B_{j}$. By \pr{e:LBac}, $B_i$ is anti-complete to each $B_j$ and so $ab\notin E(G)$. But then $G[\{v_i, a, v_{i+3},\dots, v_j,  b, v_{j+3}, \dots, v_{i-1}\}]$ is an induced  $C_{2\alpha-1}$ in $G$,  a contradiction.  \hfill\vrule height3pt width6pt depth2pt\\

With an argument similar to that  of  \pr{e:Loneneighbor}, we see that  \medskip

\noindent \refstepcounter{counter}\label{e:LNneighbor} (\arabic{counter})
For all $i\in I$, if $B_i\not=\emptyset$, then  $B_{i}$ is anti-complete to $A_j$ for any  $j\in I\less\{i-1, i, i+1, i+2\}$.\\

We next show that \medskip

\noindent \refstepcounter{counter}\label{e:Lpartition} (\arabic{counter})
For all $i\in I$, if $A_i\ne\emptyset$, then each vertex in $A_i$ is either anti-complete to $A_{i+2}$ or anti-complete to $A_{i-2}$.

\pf Suppose there exists a vertex $x\in A_i$ such that  $x$ is adjacent to a vertex $y\in A_{i-2}$ and a vertex $z\in A_{i+2}$. But then $G[ \{x, y,z\}\cup (V(C)\less \{v_{i-1}, v_i, v_{i+1}, v_{i+2}, v_{i+3}\})]$ is  an induced  $C_{2\alpha-1}$ in $G$, a contradiction. \hfill\vrule height3pt width6pt depth2pt\medskip

\noindent \refstepcounter{counter}\label{e:LA{j-1}A{j+2}} (\arabic{counter})
For all $i\in I$,    if $B_i\ne\emptyset$,  then every vertex in $B_i$ is either complete to $A_{i-1}$ or complete to  $A_{i+2}$.

\pf Suppose for a contradiction,  say $B_2\ne\emptyset$,  and  there exists a vertex $b\in B_2$  such that $b$ is not adjacent to a vertex $a_1\in A_1$ and a vertex $a_4\in A_4$. By \pr{e:LAac}, $A_1$ is anti-complete to $A_4$. Thus $G$ contains an independent set $\{b, a_1, a_4,v_0\}$ of size four when $\alpha=3$ or an independent   set $\{b, a_1, a_4, v_0, v_7, v_9, \dots, v_{2\alpha-1}\}$ of size $\alpha+1$ when $\alpha\ge4$, a contradiction. \hfill\vrule height3pt width6pt depth2pt\medskip

\noindent \refstepcounter{counter}\label{e:L3B} (\arabic{counter})
There exists an $i\in I$ such that $B_j= \emptyset$ for all $j\in I\less \{i, i+1, i+2\}$.

\pf This is obvious  if $B_k=\emptyset$ for all $k\in I$. So we may assume that $B_k\ne\emptyset$ for some $k\in I$, say $B_2\neq \emptyset$.  Then  by \pr{e:Loneneighbor},  $B_j=\emptyset$  for all  $j=5, 6, \dots, 2\alpha$.    By \pr{e:Loneneighbor} again, either $B_0\ne \emptyset$ or $B_4\ne \emptyset$ but not both.  By symmetry,  we may assume  that $B_4=\emptyset$. Similarly, either $B_0\ne \emptyset$ or $B_3\ne \emptyset$ but not both. Thus either $B_j=\emptyset$ for all $j\in I\less \{0,1,2\}$ or  $B_j=\emptyset$ for all $j\in I\less \{1,2, 3\}$. \hfill\vrule height3pt width6pt depth2pt\medskip

By \pr{e:L3B}, we may assume that  $B_j=\emptyset$ for all $j\in I\less \{1,2,3\}$. 
For any  $A_i\ne\emptyset$, where $i\in I$, 
 let $A_i^1=\{a\in A_i: a  \text{  has a neighbor in } A_{i-2}\}$, $A_i^3=\{a\in A_i: a\, \text{ has a neighbor in } A_{i+2}\}$,  and $A_i^2=A_i\less (A_i^1\cup A_i^3)$. Then $A_i^2$ is anti-complete to $A_{i-2}\cup A_{i+2}$. By \pr{e:Lpartition}, 
  $A_i^1$ is anti-complete to $A_{i+2}$ and $A_i^3$ is anti-complete to $A_{i-2}$.   Clearly, $\{A_i^1, A_i^2, A_i^3\}$ partitions $A_i$. 
  Next,  for any  $B_j\ne\emptyset$, where $j\in \{1,2,3\}$, by \pr{e:LA{j-1}A{j+2}}, let $B_j^1=\{b\in B_j: b\, \text{ is complete to }  A_{j-1}\}$ and  $B_j^2=\{b\in B_j: b\, \text{ is complete to }  A_{j+2}\}$. Clearly,  $B_j^1$ and $B_j^2$ are not necessarily disjoint. It is worth noting  that  $B_j^1$ and $B_j^2$ are not symmetrical because  $B_j^1$  is complete to $A_{j-1}$ and $B_j^2$  is complete to $A_{j+2}$.  \medskip
   
 \noindent \refstepcounter{counter}\label{e:LBjA{j-1}} (\arabic{counter})
 For any $j\in\{1,2,3\}$, $B_j$ is anti-complete to $A_{j-1}^1\cup A_{j+2}^3$.  
 
 \pf Suppose  there exist a vertex $b\in B_j$ and a vertex $a\in A_{j-1}^1\cup A_{j+2}^3$ such that $ba\in E(G)$. 
 By the definitions of $A_{j-1}^1$ and $A_{j+2}^3$, we see that $a$ has a neighbor, say $c$, in $A_{j-3}$ if $a\in A_{j-1}^1$, or  in $A_{j+4}$ if $a\in A_{j+2}^3$. But then, either  $G[ \{b, a,c\}\cup (V(C)\less \{v_{j-2}, v_{j-1}, v_j, v_{j+1}, v_{j+2}\})]$ is  an induced  $C_{2\alpha-1}$ in $G$ when    $a\in A_{j-1}^1$, or  $G[ \{b, a,c\}\cup (V(C)\less \{v_{j+1}, v_{j+2}, v_{j+3}, v_{j+4}, v_{j+5}\})]$  is  an induced  $C_{2\alpha-1}$ in $G$ when   $a\in A_{j+2}^3$. In either case, we obtain a contradiction. \hfill\vrule height3pt width6pt depth2pt\\

 \noindent \refstepcounter{counter}\label{e:Lquasi-line} (\arabic{counter})
 $G$ is a quasi-line graph.
 
 \pf   It suffices to show that for any $x\in V(G)$, $N(x)$ is covered by two  cliques.   By \pr{e:LJ}, $J=\emptyset$. Since $B_j^1$ and $B_j^2$ are not symmetrical  for all $j\in\{1,2,3\}$, we consider the following four cases. \\
 
\noindent{\bf Case 1:}   $x\in A_i$ for some $i\in I\less\{0, 1,2,3,4,5\}$.\medskip

In this case, $x\in A_i^k$ for some $k\in\{1,2,3\}$. We first assume that $k=1$. Then $x\in A_i^1$. By \pr{e:Lpartition} and the definition of $A_i^1$,  $x$ is anti-complete to $A_{i+2}$ and so  $N[x]\subseteq A_{i-2}\cup A_{i-1}\cup A_i\cup A_{i+1}\cup\{v_{i}, v_{i+1}, v_{i+2}\}$. We see that $N(x)$ is covered by two  cliques $G[(A_{i-2}\cap N(x))\cup A_{i-1}\cup \{v_i\}]$ and $G[(A_i\less x)\cup A_{i+1}\cup\{v_{i+1}, v_{i+2}\}]$. By symmetry, the same holds if $k=3$. So we may assume that $k=2$. By the definition of $A_i^2$, $x$ is anti-complete to $A_{i-2}\cup A_{i+2}$. Thus  $N[x]=A_{i-1}\cup A_i\cup A_{i+1}\cup\{v_i, v_{i+1}, v_{i+2}\}$ and so $N(x)$ is covered by two  cliques  $G[A_{i-1}\cup (A_i\less x)\cup\{v_i\}]$  and $G[A_{i+1}\cup \{v_{i+1}, v_{i+2}\}]$.   \\

\noindent{\bf Case 2:}   $x\in A_i$ for some $i\in \{0,1,2,3,4,5\}$.\medskip

In this case,  we first assume that $i=0$. Then $x\in A_0^k$ for some $k\in\{1,2,3\}$. Assume that  $x\in A_0^1$. Then  $x$ is anti-complete to $B_1$ by \pr{e:LBjA{j-1}} and anti-complete to $A_2$ by \pr{e:Lpartition}. Thus  $N[x]\subseteq A_{2\alpha-1}\cup A_{2\alpha}\cup A_0\cup A_1\cup \{v_0, v_1, v_2\}$.  One can see that $N(x)$ is covered by two  cliques   $G[(A_{2\alpha-1}\cap N(x))\cup A_{2\alpha}\cup \{v_0\}]$ and $G[(A_0\less x)\cup A_1\cup\{v_1, v_2\}]$.  It can be easily checked that $N(x)$ is covered by two  cliques $G[A_{2\alpha}\cup (A_0\less x)\cup\{v_0\}]$ and $G[A_1\cup (B_1\cap N(x))\cup\{v_1, v_2\}]$ if $k=2$; and  by two  cliques  $G[A_{2\alpha}\cup (A_0\less x)\cup\{v_0\}]$ and $G[A_1\cup (B_1\cap N(x))\cup(A_2\cap N(x))\cup\{v_1, v_2\}]$ if $k=3$.  \medskip

Next assume that $i=1$. Then $x\in A_1^k$ for some $k\in\{1,2,3\}$ and $N[x]\subseteq A_{2\alpha}\cup A_0\cup A_1\cup B_1\cup A_2\cup B_2\cup A_3\cup\{v_1, v_2, v_3\}$.  One can see that $N(x)$ is covered by two  cliques   $G[(A_{2\alpha}\cap N(x))\cup A_0\cup \{v_1\}]$ and $G[(A_1\less x)\cup B_1\cup A_2\cup\{v_2, v_3\}]$ if $k=1$; by two  cliques  $G[A_0\cup (A_1\less x)\cup \{v_1\}]$ and $G[B_1\cup A_2\cup (B_2\cap N(x))\cup\{v_2, v_3\}]$ if $k=2$; and two  cliques  $G[A_0\cup (A_1\less x)\cup B_1^1\cup\{v_1, v_2\}]$ and $G[B_1^2\cup A_2\cup (B_2\cap N(x))\cup (A_3\cap N(x))\cup\{ v_3\}]$ if $k=3$. \medskip

Assume that $i=2$. Then $x\in A_2^k$ for some $k\in\{1,2,3\}$ and $N[x]\subseteq A_{0}\cup A_{1}\cup B_{1}\cup A_2\cup B_2\cup A_{3}\cup B_{3}\cup A_4\cup\{ v_2, v_3, v_4\}$.  One can check  that $N(x)$ is covered by two  cliques   $G[(A_0\cap N(x))\cup A_1\cup B_1^1\cup \{v_2\}]$ and $G[(A_2\less x)\cup B_1^2\cup B_2\cup A_3\cup\{v_3, v_4\}]$ if $k=1$; by two  cliques  $G[A_1\cup B_1\cup (A_2\less x)\cup \{v_2\}]$ and $G[B_2\cup A_3\cup (B_3\cap N(x))\cup\{v_3, v_4\}]$ if $k=2$;  and by two cliques  $G[A_1\cup B_1\cup (A_2\less x)\cup B_2^1\cup\{v_2, v_3\}]$ and $G[B_2^2\cup A_3\cup (B_3\cap N(x))\cup (A_4\cap N(x))\cup\{ v_4\}]$ if $k=3$.  \medskip

Assume that $i=3$. Then $x\in A_3^k$ for some $k\in\{1,2,3\}$ and $N[x]\subseteq A_1\cup B_{1}\cup A_2\cup B_2\cup A_{3}\cup B_{3}\cup A_4\cup A_5\cup\{ v_3, v_4, v_5\}$.  One can check  that $N(x)$ is covered by two  cliques   $G[(A_1\cap N(x)) \cup (B_1\cap N(x))\cup A_2\cup B_2^1\cup \{v_3\}]$ and $G[B_2^2\cup (A_3\less x)\cup   B_3\cup A_4\cup\{v_4, v_5\}]$ if $k=1$; and by two  cliques  $G[(B_1\cap N(x))\cup A_2\cup B_2\cup \{v_3\}]$ and $G[(A_3\less x)\cup B_3 \cup A_4\cup\{v_4, v_5\}]$ if $k=2$. So we may assume that $x\in A_3^3$. By \pr{e:LBjA{j-1}}, we have $A_3^3$ is anti-complete to $B_1$. Thus $N(x)$ is covered by  two  cliques  $G[A_2\cup B_2\cup (A_3\less x)\cup B_3^1\cup\{v_3, v_4\}]$ and $G[B_3^2\cup A_4\cup (A_5 \cap N(x))\cup\{ v_5\}]$ if $k=3$.  \medskip

Assume that $i=4$. Then $x\in A_4^k$ for some $k\in\{1,2,3\}$ and $N[x]\subseteq A_2\cup B_{2}\cup A_3\cup B_3\cup  A_4\cup A_5\cup A_6\cup\{ v_4, v_5, v_6\}$.  One can see  that $N(x)$ is covered by two  cliques   $G[(A_2\cap N(x)) \cup (B_2\cap N(x))\cup A_3\cup B_3^1\cup \{v_4\}]$ and $G[B_3^2\cup (A_4\less x)\cup  A_5\cup\{v_5, v_6\}]$ if $k=1$; and  by two  cliques  $G[(B_2\cap N(x))\cup A_3\cup B_3\cup \{v_4\}]$ and $G[(A_4\less x) \cup A_5\cup\{v_5, v_6\}]$ if $k=2$. So we may assume that $x\in A_4^3$. By  \pr{e:LBjA{j-1}}, we have $A_4^3$ is anti-complete to $B_2$.  Thus $N(x)$ is covered  by  two  cliques  $G[A_3\cup B_3\cup (A_4\less x)\cup\{v_4, v_5\}]$ and $G[ A_5\cup (A_6 \cap N(x))\cup\{ v_6\}]$ if $k=3$. \medskip

Finally assume that $i=5$. Then $x\in A_5^k$ for some $k\in\{1,2,3\}$ and $N[x]\subseteq A_3\cup B_{3}\cup A_4\cup   A_5\cup A_6\cup A_7\cup\{ v_5, v_6, v_7\}$.  One can check  that $N(x)$ is covered by two  cliques   $G[(A_3\cap N(x)) \cup (B_3\cap N(x))\cup A_4\cup  \{v_5\}]$ and $G[(A_5\less x)\cup  A_6\cup\{v_6, v_7\}]$ if $k=1$;  and  by two  cliques  $G[(B_3\cap N(x))\cup A_4\cup \{v_5\}]$ and $G[(A_5\less x) \cup A_6\cup\{v_6, v_7\}]$ if $k=2$.  So we may assume that $x\in A_5^3$. By  \pr{e:LBjA{j-1}}, we have $A_5^3$ is anti-complete to $B_3$.   Thus $N(x)$ is covered by  two  cliques  $G[ A_4\cup (A_5\less x)\cup \{v_5, v_6\}]$ and $G[ A_6\cup (A_7\cap N(x))\cup\{v_7\}]$  if $k=3$. \medskip

This completes the proof of Case 2. \\

\noindent{\bf Case 3:}   $x\in B_j$ for some $j\in \{1,2,3\}$.\medskip

In this case, first assume that $j=1$.  Then $x\in B_1^k$ for some $k\in\{1,2\}$, and  $N[x]\subseteq A_0\cup A_1\cup B_1\cup A_2\cup B_2\cup A_3\cup\{v_1, v_2, v_3, v_4\}$. We see that  $N(x)$ is covered by two  cliques   $G[A_0\cup A_1\cup (B_1^1\less x)\cup \{v_1, v_2\}]$ and $G[B_1^2\cup A_2\cup B_2\cup (A_3\cap N(x))\cup\{v_3, v_4\}]$ if $k=1$; and by two cliques  $G[(A_0\cap N(x))\cup A_1\cup B_1^1\cup \{v_1, v_2\}]$ and $G[(B_1^2\less x)\cup A_2\cup B_2\cup A_3\cup\{v_3, v_4\}]$ if $k=2$. \medskip

Next assume that $j=2$. Then $x\in B_2^k$ for some $k\in\{1,2\}$, and  $N[x]\subseteq A_1\cup B_1\cup A_2\cup B_2\cup A_3\cup B_3\cup A_4\cup \{ v_2, v_3, v_4, v_5\}$. One can  see that $N(x)$ is covered by two cliques   $G[A_1\cup B_1\cup A_2\cup (B_2^1\less x)\cup \{v_2, v_3\}]$ and $G[B_2^2\cup A_3\cup B_3\cup (A_4\cap N(x))\cup\{v_4, v_5\}]$ if $k=1$; and by two  cliques $G[(A_1\cap N(x))\cup B_1\cup A_2\cup B_2^1\cup \{v_2, v_3\}]$ and $G[(B_2^2\less x)\cup A_3\cup B_3\cup A_4\cup\{v_4, v_5\}]$ if $k=2$. \medskip

Finally assume $j=3$. Then $x\in B_3^k$ for some $k\in\{1,2\}$, and  $N[x]\subseteq A_2\cup B_2\cup A_3\cup B_3\cup A_4\cup A_5\cup \{v_3, v_4, v_5, v_6\}$. We see that $N(x)$ is covered by two cliques   $G[A_2\cup B_2\cup A_3\cup (B_3^1\less x)\cup \{v_3, v_4\}]$ and $G[B_3^2\cup A_4\cup (A_5\cap N(x))\cup\{v_5, v_6\}]$ if $k=1$; and by two  cliques  $G[(A_2\cap N(x))\cup B_2\cup A_3\cup B_3^1\cup \{v_3, v_4\}]$ and $G[(B_3^2\less x)\cup A_4\cup A_5\cup\{v_5, v_6\}]$ if $k=2$. \\

\noindent{\bf Case 4:}   $x\in V(C)$. \medskip

In this case, let $x=v_i$ for some $i\in I$.   First assume that $i\ne 1,2,3,4,5,6$.  Then  $N(v_{i})=A_{i-2}\cup A_{i-1}\cup A_{i}\cup\{v_{i-1}, v_{i+1}\}$  and so   $N(v_{i})$ is  covered by two   cliques $G[A_{i-2}\cup A_{i-1}\cup\{v_{i-1}\}]$ and $G[A_{i}\cup\{ v_{i+1}\}]$. Next assume that $i\in\{ 1,2,3,4,5,6\}$.   One can easily check that $N(v_1)$ is covered by two  cliques $G[A_{2\alpha}\cup A_0\cup\{v_0\}]$ and $G[A_1\cup B_1\cup\{v_2\}]$;  $N(v_2)$   by  two  cliques  $G[A_{0}\cup A_1\cup\{v_1\}]$ and $G[B_1\cup A_2\cup B_2\cup\{v_3\}]$;      $N(v_3)$  by two  cliques $G[A_{1}\cup B_1\cup A_2\cup\{v_2\}]$ and $G[B_2\cup A_3\cup B_3\cup\{v_4\}]$;  $N(v_4)$  by two  cliques $G[ B_1\cup A_2\cup B_2\cup\{v_3\}]$ and $G[A_3\cup B_3\cup A_4\cup\{v_5\}]$;  $N(v_5)$  by two  cliques $G[ B_2\cup A_3\cup B_3\cup\{v_4\}]$ and $G[A_4\cup  A_5\cup\{v_6\}]$;  and $N(v_6)$  by two  cliques $G[ B_3\cup A_4\cup\{v_5\}]$ and $G[A_5\cup  A_6\cup\{v_7\}]$, respectively. 
\medskip

This completes the proof of Lemma~\ref{quasi-line}. \hfill\vrule height3pt width6pt depth2pt\\

We are now ready to prove Theorem~\ref{main}.\\

\noindent{\bf Proof of Theorem~\ref{main}}:   The statement is trivially true when $s=t=2$. So we may assume that  $t\ge3$.  Then $\chi(G) =s+t-1\ge4$. Suppose for a contradiction that $G$ is not $(s,t)$-splittable.    Let  $G$ be a counterexample with $|G|$ minimum.  
 By Theorem~\ref{spgt} and the assumption that $\chi(G)>\omega(G)$, $G$ is not perfect. Since  $G$ is $\{C_4,   C_5, C_7, \dots, C_{2\alpha(G)-1}\}$-free, by Corollary~\ref{spgt1},   we see that   $G$ must contain an induced   cycle of length   $2\alpha(G)+1$.  We   claim that  $\Delta(G)\le |G|-2$.  Suppose there exists a vertex $x$ in $G$ with $d(x)=|G|-1$. Then $\chi(G\less x)=\chi(G)-1\ge3$ because   $\chi(G)\ge4$.  Since $t\ge3$, by  minimality of $|G|$, $ G\less x$ is $(s, t-1)$-splittable. Let $S,  T'$ be  a partition of $V(G\less x)$ such that $\chi(G[S])\ge s$ and $\chi(G[T'])\ge t-1$. Let $T:=T'\cup \{x\}$. Then $S, T$ is a partition of $V(G)$ with $\chi(G[S])\ge s$ and $\chi(G[T])\ge t $.  Thus      $ G$ is $(s,t)$-splittable,  a contradiction.  This proves that  $\Delta(G)\le |G|-2$, as claimed. 
By Lemma~\ref{quasi-line}, $G$ is a quasi-line graph. By Theorem~\ref{quasiline}, $ G$ is $(s,t)$-splittable, a contradiction. This completes the proof of Theorem~\ref{main}. \hfill\vrule height3pt width6pt depth2pt\

\section*{Acknowledgement} The author  would like to thank Jason Bentley for helpful discussion.

\end{document}